# Geometric Programming Problem with Co-Efficients and Exponents Associated with Binary Numbers

A.K.Ojha[1] and A.K.Das[2]

[1] School of Basic Sciences, I.I.T Bhubaneswar
Orissa-751013, India

[2] Department of Mathematics, Bhadrak College
Bhadrak, Orissa-756100, India

**Abstract**

Geometric programming (GP) provides a power tool for solving a variety of optimization problems. In the real world, many applications of geometric programming (GP) are engineering design problems in which some of the problem parameters are estimating of actual values. This paper develops a solution procedure to solve nonlinear programming problems using GP technique by splitting the cost coefficients, constraint coefficients and exponents with the help of binary numbers. The equivalent mathematical programming problems are formulated to find their corresponding value of the objective function based on the duality theorem. The ability of calculating the cost coefficients, constraint coefficients and exponents developed in this paper might help lead to more realistic modeling efforts in engineering design areas. Standard nonlinear programming software has been used to solve the proposed optimization problem. Two numerical examples are presented to illustrate the method.

## 1. Introduction

Geometric programming (GP) a technique developed for solving algebraic non-linear programming problems subject to linear or nonlinear constraints is useful in the study of a variety of optimization problems. Dufffin, Peterson and Zener [8] put a foundation stone to solve wide range of engineering problems by developing basic theories of geometric programming and its application in their text book. Geometric programming derives its name from its intimate connection with geometrical concepts because the method is based on geometric inequality and their properties that relate sums and products of positive numbers. The application of geometric inequality has also been very useful in the construction of the condensation technique for posynomial problems [21, 26]. The notion of a condensed posynomial problem is introduced by Duffin [8] and from then, it has played an important role in the method that they use to approach the resolution to the problem from the primal perspective. One of the remarkable properties of geometric programming is that a problem with highly nonlinear constraints can be stated equivalently as one with only linear constraints. This is because there is a strong duality theorem for geometric programming problems. If the primal problem is in the posynomial form then a global minimizing solution to that problem can be obtained by solving the dual maximization. The dual constraints are linear and linearly constrained programs are generally easier to solve than ones with nonlinear constraints. GP problem has a dual impact in the area of integrated circuit design [4, 10, 20], manufacturing system design [3, 8], project management [28], maximization of long run and short run profit [18], generalized geometric programming problem with non-positive variables [34] and goal programming model [1]. Several algorithms due to Beightler and Phillips [2], Fang et al. [9], Kortanek et al. [13], Peterson [22], Rajgopal and Bricker [25] and Jhu and Kortanek [35] strengthen the solution of complicated GPP for the exact known value of the cost and constraint coefficients. Sensitive analysis of various optimal solutions due to Dembo [6], Dinkel and Tretter [7] and Kyparisis [15] using GP technique simplifies certain engineering design problem in which some of the problem parameters are estimating of actual values [2]. Entropy based problem due to Samanta and Mazumder [27] extend the application of GP to variety of entropy based transportation problem. In a recent paper Jana Mazumder [11] applied GP techniques to solve the entropy maximization problem and GP with Fuzzy parameters due to Liu [19] accelerate the work in this direction. As a few works on signomial geometric programming (SGP) have been made, various authors are concentrating to solve such problem. Some of the recent works due to Shen et al. [29] and Shen and Zharg [30] are accelerating to improve the global optimization algorithm for solving certain SGP problems. As the engineering design problems are complex in nature, the generalized GP techniques have been applied to solve such problems. Works on branch and pruning approach on GGP by Shen and Jiao [31], linearization method of global optimization due to Shen [32] and Qu et al. [24] are some of the recent applications for solving GP problems. The monomial GP with fuzzy relation due to Shivanian and Khorram [33] and maximum likelihood estimation problem proposed by Lim et al. [16] are the most recent works on GPP.

In this paper we restrict attention to posynomial geometric programs and develop a solution procedure by splitting the cost coefficients, constraint coefficients and exponents using binary numbers. This result will provide the decision makers with more information for making better decisions.

The organization of this paper is as follows: following introduction, mathematical formulations and methodology for solving GPP with binary coefficients have been discussed to find the objective value of the problem in section-2. Transformation of the problem into binary model has been discussed in section-3. Dual form of GPP has been discussed in section-4. Some illustrative examples are given in section-5 for understanding the problems and finally in section-6 some conclusions are drawn from the discussion.

## 2. Mathematical Formulation

A typical constrained posynomial geometric programming problem is presented as follows:

$$Z = \min_{x} : f_0(x) = \sum_{t=1}^{T_0} C_{0t} \prod_{j=1}^{n} x_j^{a_{0tj}}$$

subject to

$$f_i(x) = \sum_{t=1}^{T_i} C_{it} \prod_{j=1}^{n} x_j^{a_{itj}} \leq 1,$$

$$i = 1, 2, \ldots, m \quad (2.1)$$

$$x_j > 0, \; j = 1, 2, 3, \ldots, n.$$

The posynomial $f_0(x)$ is an objective function containing $T_0$ number of terms where as the posynomial $f_i(x)$, $i = 1, 2, .., m$ contains $T_i$ terms with m inequality constraints. By the definition of posynomial all the co-efficients $C_{it}$, $i = 0, 1, 2, \ldots, m$ and $t = 1, 2, \ldots, T_i$ are positive and the exponents $a_{0tj}$ and $a_{itj}$ are arbitrary constants. Writing the right hand side of the geometric programming problem given by (2.1) in more general form, we have

$$Z = \min_{x} : f_0(x) = \sum_{t=1}^{T_0} C_{0t} \prod_{j=1}^{n} x_j^{a_{0tj}}$$

subject to

$$f_i(x) = \sum_{t=1}^{T_i} C_{it} \prod_{j=1}^{n} x_j^{a_{itj}} \leq b_i$$

$$i = 1, 2, \ldots, m \quad (2.2)$$

$$x_j > 0, \; j = 1, 2, \ldots n.$$

where all $b_i$ are positive real numbers. If $b_i = 1$ for all i then this modified geometric program becomes the original one given by (2.1).

Since $b_i$ in the model (2.2) may not be equal to constant 1 then dividing the constant coefficients $C_{it}$ by $b_i \; \forall \; i$ then it is transformed to the standard form

$$Z = \min_{x} : f_0(x) = \sum_{t=1}^{T_0} C_{0t} \prod_{j=1}^{n} x_j^{a_{0tj}}$$

subject to

$$f_i(x) = \sum_{t=1}^{T_i} C_{it}(b_i)^{-1} \prod_{j=1}^{n} x_j^{a_{itj}} \leq 1$$

$$i = 1, 2, \ldots, m \quad (2.3)$$

$$x_j > 0, j = 1, 2, \ldots n.$$

In this paper our aim is to solve the above problem by splitting the coefficients and exponents of the variables using binary numbers. In order to solve this problem it is necessary to transform the problem to the standard mathematical programming problem using corresponding binary expressions. The following section explains the detail procedure for converting the problem using binary numbers in the different uses.

## 3. Transformation of the Problem into Binary Model

The proposed GP given by (2.1) and (2.2) can be reformulated by splitting any of $C_{ot}$, $C_{it}$, and $a_{itj}$ to the 3 to 8 terms using binary numbers.

**Case– 1:** When $C_{ot}$ has 3 terms.

$$Z = \min_{x} : \sum_{t=1}^{T_0}(a_1^t z_1(1 - z_2) + a_2^t (1 - z_1)z_2 + a_3^t (1-z_1)(1 - z_2)) \prod_{j=1}^{n} x_j^{a_{otj}}$$

subject to

$$\sum_{t=1}^{T_i}(a_1^t z_1(1 - z_2) + a_2^t (1 - z_1)z_2 + a_3^t (1-z_1)(1 - z_2)) \prod_{j=1}^{n} x_j^{a_{itj}} \leq 1$$

$$i = 1, 2, \ldots, m \quad (3.1)$$
$$z_1 + z_2 \leq 1 \quad (3.2)$$

**Case -2:** When $C_{ot}$ has 4 terms

$$\min_{x} : \quad T_0$$

51

$$Z = \sum_{t=1} (a_1^t z_1(1-z_2) + a_2^t (1-z_1)z_2 + a_3^t (1-z_1)(1-z_2) + a_4^t z_1 z_2) \prod_{j=1}^n x_j^{a_{otj}}$$

subject to

$$\sum_{t=1}^{T_i} (a_1^t z_1(1-z_2) + a_2^t (1-z_1)z_2 + a_3^t (1-z_1)(1-z_2) + a_4^t z_1 z_2) \prod_{j=1}^n x_j^{a_{itj}} \leq 1, \; i = 1, 2, \ldots, m \quad (3.3)$$

$$z_1 + z_2 \leq 1 \quad (3.4)$$

**Case -3 :** When $C_{0t}$ has 5 terms.

$$Z = \min_x : \sum_{t=1}^{T_0} (a_1^t z_1(1-z_2)(1-z_3) + a_2^t (1-z_1)z_2(1-z_3) + a_3^t (1-z_1)(1-z_2)z_3 + a_4^t (1-z_1)(1-z_2)(1-z_3) + a_5^t z_1 z_2 z_3) \prod_{j=1}^n x_j^{a_{otj}}$$

subject to

$$\sum_{t=1}^{T_i} (a_1^t z_1(1-z_2)(1-z_3) + a_2^t (1-z_1)z_2(1-z_3) + a_3^t (1-z_1)(1-z_2)z_3 + a_4^t (1-z_1)(1-z_2)(1-z_3) + a_5^t z_1 z_2 z_3) \prod_{j=1}^n x_j^{a_{itj}} \leq 1,$$

$i = 1, 2, \ldots m$ (3.5)
$z_1 z_2 (1-z_3) = 0$ (3.6)
$z_2 z_3 (1-z_1) = 0$ (3.7)
$z_1 z_3 (1-z_2) = 0$ (3.8)

**Case -4 :** When $C_{ot}$ has 6 terms.

$$Z = \min_x : \sum_{t=1}^{T_0} (a_1^t z_1(1-z_2)(1-z_3) + a_2^t (1-z_1)z_2(1-z_3) + a_3^t (1-z_1)(1-z_2)z_3 + a_4^t z_1 z_2 (1-z_3) + a_5^t z_1 z_3(1-z_2) + a_6^t z_2 z_3(1-z_1)) \prod_{j=1}^n x_j^{a_{otj}}$$

subject to

$$\sum_{t=1}^{T_i} (a_1^t z_1(1-z_2)(1-z_3) + a_2^t (1-z_1)z_2(1-z_3) + a_3^t (1-z_1)(1-z_2)z_3 + a_4^t z_1 z_2 (1-z_3) + a_5^t z_1 z_3(1-z_2) + a_6^t z_2 z_3(1-z_1)) \prod_{j=1}^n x_j^{a_{itj}} \leq 1, \; i = 1, 2, \ldots m \quad (3.9)$$

$z_1 z_2 z_3 = 0$ (3.10)
$(1-z_1)(1-z_2)(1-z_3) = 0$ (3.11)

**Case -5 :** When $C_{0t}$ has 7 terms

$$Z = \min_x : \sum_{t=1}^{T_0} (a_1^t z_1(1-z_2)(1-z_3) + a_2^t (1-z_1)z_2(1-z_3) + a_3^t (1-z_1)(1-z_2)z_3 + a_4^t z_1 z_2 (1-z_3) + a_5^t z_1 z_3(1-z_2) + a_6^t z_2 z_3(1-z_1) + a_7^t (1-z_1)(1-z_2)(1-z_3)) \prod_{j=1}^n x_j^{a_{otj}}$$

subject to

$$\sum_{t=1}^{T_i} (a_1^t z_1(1-z_2)(1-z_3) + a_2^t (1-z_1)z_2(1-z_3) + a_3^t (1-z_1)(1-z_2)z_3 + a_4^t z_1 z_2 (1-z_3) + a_5^t z_1 z_3(1-z_2) + a_6^t z_2 z_3(1-z_1) + a_7^t (1-z_1)(1-z_2)(1-z_3)) \prod_{j=1}^n x_j^{a_{itj}} \leq 1,$$

$i = 1, 2, \ldots, m$ (3.12)
$z_1 z_2 z_3 = 0$ (3.13)

**Case -6:** When $C_{ot}$ has 8 terms.

$$Z = \min_x : \sum_{t=1}^{T_0} (a_1^t z_1(1-z_2)(1-z_3) + a_2^t (1-z_1)z_2(1-z_3) + a_3^t (1-z_1)(1-z_2)z_3 + a_4^t z_1 z_2 (1-z_3) + a_5^t z_1 z_3(1-z_2) + a_6^t z_2 z_3(1-z_1) + a_7^t (1-z_1)(1-z_2)(1-z_3) + a_8^t z_1 z_2 z_3) \prod_{j=1}^n x_j^{a_{otj}}$$

subject to :

$$\sum_{t=1}^{T_i} (a_1^t z_1(1-z_2)(1-z_3) + a_2^t (1-z_1)z_2(1-z_3) + a_3^t (1-z_1)(1-z_2)z_3 + a_4^t z_1 z_2 (1-z_3)$$

IJCSI

$$+ a_3^t (1-z_1)(1-z_2) z_3 + a_4^t z_1 z_2 (1-z_3)$$

$$+ a_5^t z_1 z_3 (1-z_2) + a_6^t z_2 z_3 (1-z_1)$$

$$+ a_7^t (1-z_1)(1-z_2)(1-z_3)$$

$$+ a_8^t z_1 z_2 z_3) \prod_{j=1}^{n} x_j^{a_{otj}} \leq 1,$$

$i = 1, 2, \ldots m$ (3.15)

$z_1 z_2 z_3 = 0$ (3.16)

## 4. Dual form of GPP

Since model (2.2) is the conventional geometric programming problem then it can be solved directly by using primal based algorithm or dual based algorithm [22]. Methods due to Rajgopal and Bricker [25], Beightler and Phillips [2] and Duffin [8] projected in their analysis that the dual problem has the desirable features of being linearly constrained and having an objective function with attractive structural properties with more suitable solution. According to Duffin [8], Beightler and Phillips [2] one can transform the program of (2.2) into the corresponding dual geometric problem as follows:

$$Z = \max_{w} : \prod_{t=1}^{T_0} \left(\frac{C_{0t}}{w_{0t}}\right)^{w_{0t}} \prod_{i=1}^{m} \prod_{t=1}^{T_i} \left(\frac{C_{it} b_i^{-1} w_{i0}}{w_{it}}\right)^{w_{it}}$$

$$\prod_{t=1}^{T_i} \lambda(w_{it})^{\lambda(w_{it})}$$

$$\sum_{t=1}^{T_0} w_{0t} = 1 \quad (4.1)$$

$$\sum_{i=1}^{m} \sum_{t=1}^{T_i} a_{itj} w_{it} = 0, \quad j = 1, 2, \ldots n. \quad (4.2)$$

$w_{it} \geq 0 \quad \forall t, i.$ (4.3)

The model (4.1) is the useful dual problem and it can be solved using the method relating to the dual theorem.

## 5. Illustrative Examples

To illustrate the methodology proposed in this paper for solving a GPP when splitting the cost coefficients, constraint coefficients and the exponents of the decision variables using binary numbers, a few numerical examples are considered.

**Example: 1**

Let us consider the geometric programming problem which has the following mathematical form:

$$z = \min_{x} : cx_1^p + \frac{3}{x_2^3} + x_1 x_2 \quad (5.1)$$

subject to $ax_1 + x_2 \leq 1$

$x_1, x_2 > 0$

**Primal Solution**

**Case -1:** Splitting c, p, a into 3 terms we have

min = c * $x_1$ ∧ p + 3 / ($x_2$∧3)+$x_1$* $x_2$; (5.2)

a * $x_1$ + $x_2$ <= 1;

c = 5 * (1-$z_1$) * (1-$z_2$) + (1-$z_1$) * $z_2$
+ 3 * $z_1$ * (1-$z_2$) (5.3)

p=-(1-$z_3$)*(1-$z_4$)–3*(1-$z_3$)*$z_4$-2*$z_3$*(1-$z_4$) (5.4)

a=3*(1-$z_5$)*(1-$z_6$)+(1-$z_5$)*$z_6$+2*$z_5$*(1-$z_6$) (5.5)

$z_1 + z_2$ 1

$z_3 + z_4$ 1

$z_5 + z_6$ 1

**Case–2:** Splitting c, p, a into 4 terms we have

min = c * $x_1$ ∧ p + 3 / ($x_2$ ∧ 3) + $x_1$ * $x_2$ (5.6)

a * $x_1$ + $x_2$ 1

c=5*(1-$z_1$)*(1-$z_2$)+(1-$z_1$)*$z_2$+3*$z_1$*(1-$z_2$)+4*$z_1$*$z_2$ (5.7)

p=-(1-$z_3$)*(1-$z_4$)–3*(1-$z_3$)*$z_4$-2*$z_3$*(1-$z_4$)-4*$z_3$*$z_4$ (5.8)

a=3*(1-$z_5$)*(1-$z_6$)+(1-$z_5$)*$z_6$+2*$z_5$*(1-$z_6$)+4*$z_5$*$z_6$ (5.9)

$z_1 + z_2$ <= 1

$z_3 + z_4$ <= 1

$z_5 + z_6$ <= 14

**Case- 3 :** Splitting c, p, a into 5 terms we have

min = c * $x_1$ ∧ p + 3 / ($x_2$ ∧ 3) + $x_1$ * $x_2$ (5.10)

a * $x_1$ + $x_2$ <= 1

c=5*$z_1$*(1-$z_2$)*(1-$z_3$)+$z_2$*(1-$z_1$)*(1-$z_3$)+3*$z_3$*(1-$z_1$)*(1-$z_2$)+4*(1-$z_1$)*(1-$z_2$)*(1-$z_3$)+6*$z_1$*$z_2$ $z_3$ (5.11)

p=-$z_4$*(1-$z_5$)*(1-$z_6$)–3*$z_5$*(1-$z_4$)*(1-$z_6$)-2*$z_6$*(1-$z_4$)*(1-$z_5$)–(1-$z_4$)*(1-$z_5$)*(1-$z_6$)-4*$z_4$*$z_5$*$z_6$ (5.12)

a=3*$z_7$*(1-$z_8$)*(1-$z_9$)+$z_8$*(1-$z_7$)*(1-$z_9$)+2*$z_9$*(1-$z_7$)*(1-$z_8$)+5*(1-$z_7$)*(1-$z_8$)*(1-$z_9$)+6*$z_7$*$z_8$*$z_9$ (5.13)

$z_1$ * $z_2$ * (1-$z_3$) = 0

$z_2$ * $z_3$ * (1-$z_4$) = 0

$z_1$ * $z_3$ * (1-$z_2$) = 0

$z_4$*$z_5$*(1-$z_6$)=0

$z_5$ * $z_6$ * (1-$z_4$) = 0;

$z_4$ * $z_6$ * (1-$z_5$) = 0;

$z_7$ * $z_8$ * (1-$z_9$) = 0;

$z_8$ * $z_9$ * (1-$z_7$) = 0;

$z_7$ * $z_9$ * (1-$z_8$) = 0;

**Case -4:** Splitting c, p, a into 6 terms we have

min =c*$x_1$∧p+3/($x_2$ ∧ 3)+$x_1$*$x_2$; (5.14)

a * $x_1$ + $x_2$ <= 1;

$c=5*z_1*(1-z_2)*(1-z_3)+z_2*(1-z_1)*(1-z_3)+3*z_3*(1-z_1)*(1-z_2)+4z_1*z_2*(1-z_3)+6z_1*z_3*(1-z_2)+2*z_2*z_3*(1-z_1);$ (5.15)

$p=-z_4*(1-z_5)*(1-z_6)–3*z_5*(1-z_4)*(1-z_6)-2*z_6*(1-z_4)*(1-z_5)–z_4*z_5*(1-z_6)–4z_4*z_6*(1-z_5)-5*z_5*z_6*(1-z_4);$ (5.16)

$a=3*z_7*(1-z_8)*(1-z_9)+z_8*(1-z_7)*(1-z_9)+2*z_9*(1-z_7)*(1-z_8)+5*z_7*z_8*(1-z_9)+6*z_7*z_9*(1-z_8)+4*z_8*z_9*(1-z_7);$ (5.17)

$z_1 * z_2 * z_3 = 0;$

$(1-z_1) * (1-z_2) * (1-z_3) = 0;$

$z_4 * z_5 * z_6 = 0;$

$(1-z_4) * (1-z_5) * (1-z_6) = 0;$

$z_7 * z_8 * z_9 = 0;$

$(1-z_7) * (1-z_8) * (1-z_9) = 0;$

**Case-5:** Splitting c, p, a into 7 terms we have

min=$c*x_1 \wedge p+3/(x_2\wedge 3)+x_1*x_2;$ (5.18)

$a*x_1+x_2$ 1;

$c=5*z_1*(1-z_2)*(1-z_3)+z_2*(1-z_1)*(1-z_3)+3*z_3*(1-z_1)*(1-z_2)+4z_1*z_2*(1-z_3)+6z_1*z_3*(1-z_2)+2*z_2*z_3*(1-z_1)+(1-z_1)*(1-z_2)*(1-z_3);$ (5.19)

$p=-z_4*(1-z_5)*(1-z_6)–3*z_5*(1-z_4)*(1-z_6)-2*z_6*(1-z_4)*(1-z_5)–z_4*z_5*(1-z_6)–4z_4*z_6*(1-z_5)-5*z_5*z_6*(1-z_4)-(1-z_4)*(1-z_5)*(1-z_6);$ (5.20)

$a=3*z_7*(1-z_8)*(1-z_9)+z_8*(1-z_7)*(1-z_9)+2*z_9*(1-z_7)*(1-z_8)+5*z_7*z_8*(1-z_9)+6*z_7*z_9*(1-z_8)+4*z_8*z_9*(1-z_7)+4*(1-z_7)*(1-z_8)*(1-z_9);$ (5.21)

$z_1 * z_2 * z_3 = 0;$

$z_4 * z_5 * z_6 = 0;$

$z_7 * z_8 * z_9 = 0;$

**Case-6:** Splitting c, p, a into 8 terms we have

min=$c*x_1\wedge p+3/(x_2\wedge 3)+ x_1*x_2;$ (5.22)

$a * x_1 + x_2 <= 1;$

$c=5*z_1*(1-z_2)*(1-z_3)+z_2*(1-z_1)*(1-z_3)+3*z_3*(1-z_1)*(1-z_2)+4z_1*z_2*(1-z_3)+6z_1*z_3*(1-z_2)+2*z_2*z_3*(1-z_1)+(1-z_1)*(1-z_2)*(1-z_3)+2*z_1*z_2*z_3;$ (5.23)

$p=-z_4*(1-z_5)*(1-z_6)–3*z_5*(1-z_4)*(1-z_6)-2*z_6*(1-z_4)*(1-z_5)–z_4*z_5*(1-z_6)–4z_4*z_6*(1-z_5)-5*z_5*z_6*(1-z_4)-(1-z_4)*(1-z_5)*(1-z_6)–2*z_4*z_5*z_6 ;$ (5.24)

$a=3*z_7*(1-z_8)*(1-z_9)+z_8*(1-z_7)*(1-z_9)+2*z_9*(1-z_7)*(1-z_8)+5*z_7*z_8*(1-z_9)+6*z_7*z_9*(1-z_8)+4*z_8*z_9*(1-z_7)+4*(1-z_7)*(1-z_8)*(1-z_9)+2*z_7*z_8*z_9;$ (5.25)

$z_1*z_2*z_3=0;$

$z_4*z_5*z_6=0;$

$z_7*z_8*z_9=0;$

Using Lingo Software we have found that the solution in each case gives z=11.01098,c=1,p=-1,a=1, $x_1$= 0.2069792, $x_2$ = 0.7930208.

**Dual Solution**

max =$(1/w_{01})\wedge w_{01}*(3/w_{02})\wedge w_{02}*(1/w_{03})\wedge w_{03}*(1/w_{11})\wedge w_{11}*(1/ w_{12})\wedge w_{12}*(w_{11}+w_{12})\wedge (w_{11}+w_{12});$ (5.26)

$w_{01} + w_{02} + w_{03} = 1;$

$-w_{01} + w_{03} + w_{11} = 0;$

$-3 * w_{02} + w_{03} + w_{12} = 0;$

$w_{01}, w_{02}, w_{03}, w_{11}, w_{12} > 0;$

In this case we shall find z = 11.01098, $w_{01}$ = 0.4387805, $w_{02}$ = 0.5463127, $w_{03}$=0.1490681E – 01, $w_{11}$ = 0.4238737, $w_{12}$ = 1.624031.

**Example: 2**

Let us consider the geometric programming problem which has the following mathematical form:

Z = $\min_{x}$ : $cx_1+10x_2+4x_3+2x_4$ (5.27)

subject to $ax_1^p x_4^{-2} + x_2^2 x_4^{-2}$ 1

$100x_1^{-1} x_2^{-1} x_3^{-1}$ 1

$x_1, x_2, x_3, x_4 > 0$

**Primal Solution**

**Case-1:** Splitting c, p, a into 3 terms we have

min=$c*x_1+10*x_2+4*x_3+2*x_4;$ (5.28)

$a*x_1\wedge p*x_4\wedge -2+x_2\wedge 2*x_4\wedge -2$ 1;

$100*x_1\wedge -1*x_2\wedge -1+x_3\wedge -1$ 1;

$c=5*(1-z_1)*(1-z_2)+(1-z_1)*z_2+3*(1-z_2)*z_1;$ (5.29)

$p=-(1-z_3)*(1-z_4)–3*(1-z_3)*z_4 -2*(1-z_4)*z_3;$ (5.30)

$a=3*(1-z_5)*(1-z_6)+(1-z_5)*z_6 +2*(1-z_6)*z_5;$ (5.31)

$z_1 + z_2 <= 1;$

$z_3 + z_4 <= 1;$

$z_5 + z_6 <= 1;$

**Case-2:** Splitting c, p, a into 4 terms we have

min=$c*x_1+10*x_2+4*x_3+2*x_4;$ (5.32)

$a*x_1\wedge p*x_4\wedge -2+x_2\wedge 2*x_4\wedge -2$ 1;

$100*x_1\wedge -1*x_2\wedge -1*x_3\wedge -1$ 1;

$c=5*(1-z_1)*(1-z_2)+(1-z_1)*z_2+3*(1-z_2)*z_1+4*z_1*z_2;$ (5.33)

$p=-(1-z_3)*(1-z_4)–3*(1-z_3)*z_4-2*(1-z_4)*z_3-4*z_3*z_4;$ (5.34)

$a=3*(1-z_5)*(1-z_6)+(1-z_5)*z_6+2*(1-z_6)*z_5+4*z_5*z_6;$ (5.35)

$z_1 + z_2 <= 1;$

$z_3 + z_4 <= 1;$

$z_5 + z_6 <= 1;$

**Case-3:** Splitting c, p, a into 5 terms we have

min=$c*x_1+10*x_2 + 4 * x_3 + 2 * x_4;$ (5.36)

$a*x_1\wedge p*x_4\wedge -2+x_2\wedge 2*x_4\wedge -2$ 1;

$100*x_1\wedge -1*x_2\wedge -1*x_3\wedge -1$ 1;

$c=5*z_1*(1-z_2)*(1-z_3)+z_2*(1-z_1)*(1-z_3)+3*z_3*(1-z_1)*(1-z_2)+3*(1-z_1)*(1-z_2)*(1-z_3)+5*z_1*z_2*z_3;$ (5.37)

$p=z_4*(1-z_5)*(1-z_6)–3*z_5*(1-z_4)*(1-z_6)-2*z_6*(1-z_4)*(1-z_5)-3*(1-z_4)*(1-z_5)*(1-z_6)–2*z_4*z_5*z_6;$ (5.38)

$a=z_7*(1-z_8)*(1-z_9)+z_8*(1-z_7)*(1-z_9)+2*z_9*(1-z_7)*(1-z_8)+4*(1-z_7)*(1-z_8)*(1-z_9)+3*z_7*z_8*z_9 ;$ (5.39)

$z_1 * z_2 * (1-z_3) = 0;$

$z_2 * z_3 * (1-z_4) = 0;$

$z_1 * z_3 * (1-z_2) = 0;$

$z_4 * z_5 * (1-z_6) = 0;$

$z_5 * z_6 * (1-z_4) = 0;$

$z_4 * z_6 * (1-z_5) = 0;$

$z_7 * z_8 * (1-z_9) = 0;$
$z_8 * z_9 * (1-z_7) = 0;$
$z_7 * z_9 * (1-z_8) = 0;$

**Case-4:** Splitting c, p, a into 6 terms we have

$\min = c*x_1+10*x_2+4*x_3+2*x_4;$ (5.40)

$a*x_1 \wedge p*x_4 \wedge -2 + x_2 \wedge 2*x_4 \wedge -2 \leq 1;$

$100*x_1 \wedge -1*x_2 \wedge -1*x_3 \wedge -1 \leq 1;$

$c=5*z_1*(1-z_2)*(1-z_3)+z_2*(1-z_1)*(1-z_3)+3*z_3*(1-z_1)*(1-z_2)+3*z_1*z_2*(1-z_3)+z_1*z_3*(1-z_2)+z_2*z_3*(1-z_1);$ (5.41)

$p=z_4*(1-z_5)*(1-z_6)-3*z_5*(1-z_4)*(1-z_6)-2*z_6*(1-z_4)*(1-z_5)-2*z_4*z_5*(1-z_6)-3*z_4*z_6(1-z_5)-z_5*z_6*(1-z_4);$ (5.42)

$a= z_7*(1-z_8)*(1-z_9)+z_8*(1-z_7)*(1-z_9)+2*z_9*(1-z_7)*(1-z_8)+4*z_7*z_8*(1-z_9)+3*z_7*z_9*(1-z_8)+4*z_8*z_9*(1-z_7);$ (5.43)

$z_1*z_2*z_3=0;$
$(1-z_1)*(1-z_2)*(1-z_3)=0;$
$z_4*z_5*z_6=0;$
$(1-z_4)*(1-z_5)*(1-z_6)=0;$
$z_7*z_8*z_9=0;$
$(1-z_7)*(1-z_8)*(1-z_9)=0;$

**Case-5:** Splitting c, p, a into 7 terms we have

$\min = c*x_1+10*x_2+4*x_3+2*x_4;$ (5.44)

$a*x_1 \wedge p*x_4 \wedge -2 + x_2 \wedge 2*x_4 \wedge -2 <= 1;$

$100*x_1 \wedge -1*x_2 \wedge -1*x_3 \wedge -1 <= 1;$

$c=5*z_1*(1-z_2)*(1-z_3)+z_2*(1-z_1)*(1-z_3)+3*z_3*(1-z_1)*(1-z_2)+3*z_1*z_2*(1-z_3)+5*z_1*z_3*(1-z_2)+2*z_2*z_3*(1-z_1)+(1-z_1)+(1-z_1)*(1-z_2)*(1-z_3);$ (5.45)

$p=-z_4*(1-z_5)*(1-z_6)-3*z_5*(1-z_4)*(1-z_6)-2*z_6*(1-z_4)*(1-z_5)-3*z_4*z_5*(1-z_6)-2*z_4*z_6*(1-z_5)-z_5*z_6*(1-z_4)-(1-z_4)*(1-z_5)*(1-z_6);$ (5.46)

$a=z_7*(1-z_8)*(1-z_9)+z_8*(1-z_7)*(1-z_9)+2*z_9*(1-z_7)*(1-z_8)+4*z_7*z_8*(1-z_9)+3*z_7*z_9*(1-z_8)+4*z_8*z_9*(1-z_7)+4*(1-z_7)*(1-z_8)*(1-z_9);$ (5.47)

$z_1 * z_2 * z_3 = 0;$
$z_4 * z_5 * z_6 = 0;$
$z_7 * z_8 * z_9 = 0;$

**Case-6:** Splitting c, p, a into 8 terms we have

$\min = c*x_1 + 10*x_2 + 4*x_3 + 2*x_4;$ (5.48)

$a*x_1 \wedge p*x_4 \wedge -2 + x_2 \wedge 2*x_4 \wedge -2 <= 1;$

$100*x_1 \wedge -1*x_2 \wedge -1*x_3 \wedge -1 <= 1;$

$c=5*z_1*(1-z_2)*(1-z_3)+z_2*(1-z_1)*(1-z_3)+3*z_3*(1-z_1)*(1-z_2)+3*z_1*z_2*(1-z_3)+5*z_1*z_3*(1-z_2)+2*z_2*z_3*(1-z_1)+(1-z_1)*(1-z_2)*(1-z_3)+z_1*z_2*z_3;$ (5.49)

$p=-z_4*(1-z_5)*(1-z_6)-3*z_5*(1-z_4)*(1-z_6)-2*z_6*(1-z_4)*(1-z_5)-3*z_4*z_5*(1-z_6)-2*z_4*z_6*(1-z_5)-z_5*z_6*(1-z_4)-(1-z_4)*(1-z_5)*(1-z_6)-2*z_4*z_5*z_6;$ (5.50)

$a=z_7*(1-z_8)*(1-z_9)+z_8*(1-z_7)*(1-z_9)+2*z_9*(1-z_7)*(1-z_8)+4*z_7*z_8*(1-z_9)+3*z_7*z_9*(1-z_8)+4*z_8*z_9*(1-z_7)+4*(1-z_7)*(1-z_8)*(1-z_9)+2*z_7*z_8*z_9;$ (5.51)

$z_1 * z_2 * z_3 = 0;$
$z_4 * z_5 * z_6 = 0;$
$z_7 * z_8 * z_9 = 0;$

Using Lingo Software we have found that the solution in each case gives $z = 50.60611$, $c = 1$, $a = 1$, $p = -3$, $x_1 = 16.86890$, $x_2 = 1.405717$, $x_3 = 4.217114$, $x_4 = 1.405791$.

**Dual Solution:**

$\max = (1/w_{01}) \wedge w_{01} * (10/w_{02}) \wedge w_{02} * (4/w_{03}) \wedge w_{03} * (2/w_{04}) \wedge w_{04} * (1/ w_{11}) \wedge w_{11} * (1/w_{12}) \wedge w_{12} * (w_{11}+w_{12}) \wedge (w_{11}+w_{12}) * 100 \wedge w_{21};$ (5.52)

$w_{01} + w_{02} + w_{03} + w_{04} = 1;$
$w_{01} + (-3) * w_{11} - w_{21} = 0;$
$w_{02} + 2 * w_{12} - w_{21} = 0;$
$w_{03} - w_{21} = 0;$
$w_{04} - 2 * w_{11} - 2 * w_{12} = 0;$
$w_{01}, w_{02}, w_{03}, w_{04}, w_{11}, w_{12}, w_{21} > 0;$

In this case we find $z = 50.60611$, $w_{01} = 0.3333372$, $w_{02} = 0.2777762$, $w_{03} = 0.3333285$, $w_{04} = 0.5555811E - 01$, $w_{11} = 0.2903339E-05$, $w_{12} = 0.2777615E-01$, $w_{21} = 0.3333285$.

## 6. Conclusions

Since 1960 geometric programming problem has undergone several changes. In most of the engineering problems the parameters are considered as deterministic. In this paper we have discussed the problems by splitting the cost coefficients, constraint coefficients and exponents using binary numbers. Geometric programming has already shown its power in practice in the past. In many real world geometric programming problem the parameters may not be known precisely due to insufficient information and hence this paper will help the wider applications in the field of engineering problems.

**Dr.A.K.Ojha**: Dr A.K.Ojha received a Ph.D (Mathematics) from Utkal University in 1997.Currently he is an Asst. Prof. in Mathematics at I.I.T Bhubaneswar, India. He is performing research in Neural Network, Genetical Algorithm, Geometric Programming and Particle Swarm Optimization. He is served more than more than 27 years in different Govt. colleges in the state of Orissa. He is published 22 research paper in different journals and 7 books for degree students such as: Fortran 77 Programming, A text book of Modern Algebra, Fundamentals of Numerical Analysis etc.

**A.K.Das**: Mr.A.K.Das received a M.Sc (Mathematics) from Revenshaw College. Currently he is an Lecturer in Mathematics at Bhadrak College, Bhadrak, Orissa, India. He is performing research in Geometric Programming and Optimization Theory.